\documentclass[11pt]{article}
\usepackage{amsmath,amssymb,euscript,graphicx,amsfonts}
\usepackage{color}
\usepackage[dvips]{epsfig}
\usepackage{enumerate}

\title{HAMILTONICITY OF VERTEX-TRANSITIVE GRAPHS OF ORDER $4p$}
\date{\today}

\newtheorem{theorem}{Theorem}[section]
\newtheorem{proposition}[theorem]{Proposition}
\newtheorem{lemma}[theorem]{Lemma}

\newcommand{\Qed}{\rule{2.5mm}{3mm}}

\newcommand{\Aut}{\hbox{{\rm Aut}}}

\newcommand{\LCF}{\hbox{{\rm LCF}}}

\newcommand{\la}{\langle}
\newcommand{\ra}{\rangle}

\newcommand{\ZZ}{\mathbb{Z}}

\newcommand{\Cay}{\hbox{Cay}}

\newcommand{\PP}{{\cal P}}

\newcommand{\B}{{\cal{B}}}

\newcommand{\W}{{\cal W}}

\newcommand{\PSL}{\hbox{{\rm PSL}}}

\newcommand{\proofT}{ {\em\bf \hspace{-1mm}Proof of Theorem~\ref{the:main}:}}

\newenvironment{proof}{{\noindent \sc Proof.}}{\hfill $\Qed$ \\}

 \newcounter{case}
 \newenvironment{case}[1][\unskip]{\refstepcounter{case}\sc 
 \medskip \indent Case \thecase\ #1.\ }{\unskip\upshape}
 \renewcommand{\thecase}{\arabic{case}}

 \newcounter{subcase}
 \newenvironment{subcase}[1][\unskip]{\refstepcounter{subcase}\sc 
 \medskip \indent Subcase \thesubcase\ #1.\ }{\unskip\upshape}
\numberwithin{subcase}{case}

\def\ZZ{{\hbox{\sf Z\kern-.43emZ}}}

\input amssym.def

\input amssym.tex

\begin{document}


\begin{center}
{\bf\large HAMILTONICITY OF VERTEX-TRANSITIVE GRAPHS OF ORDER $4p$}
\end{center}
\bigskip\noindent
\begin{center}

  Klavdija Kutnar{\small$^{a,}$}\footnotemark \, and
 Dragan Maru\v si\v c{\small$^a$}$^,${\small$^b$}$^,$\footnotemark$^,$*\\

\bigskip

{\it {\small  $^a$University of Primorska, Cankarjeva 6, 6000 Koper, Slovenia\\
$^b$IMFM, University of Ljubljana, Jadranska 19, 1000 Ljubljana, Slovenia}}
\
\end{center}

\addtocounter{footnote}{-1}
\footnotetext{Supported in part by
``Agencija za raziskovalno dejavnost Republike Slovenije'', proj.  mladi raziskovalci.}
\addtocounter{footnote}{1}  \footnotetext{Supported in part by
``Agencija za raziskovalno dejavnost Republike Slovenije'', research program P1-0285.

~*Corresponding author e-mail: ~dragan.marusic@guest.arnes.si}

\begin{abstract}
It is shown that every connected vertex-transitive graph of order $4p$,
where $p$ is a prime, is hamiltonian with the exception of the Coxeter graph
which is known
to possess a Hamilton path.
\end{abstract}

\bigskip
\begin{quotation}
\end{quotation}


\bigskip
\section{Introductory remarks}
\label{sec:intro}
\indent

In 1969, Lov\'asz \cite{LL70} asked if   every
finite, connected vertex-transitive graph has a Hamilton path, that
is, a path going through all vertices of the graph.
With the exception of $K_2$, only four   connected
vertex-transitive graphs that do not have a Hamilton cycle are known to exist.
These four graphs are the Petersen graph, the Coxeter graph and the two graphs obtained
from them by replacing each vertex by a triangle.
The fact that none of these four graphs
is a Cayley graph  has led to a folklore conjecture that
every Cayley graph is hamiltonian (see \cite{BA89,
AZ89, ALW90, AQ01, CG96, DGMW98, ED83, GY96, KW85, DM83,
DW82, DW85,DW86} for the current status of this conjecture).

Coming back to vertex-transitive graphs,
it was shown in \cite{DGMW98} that, with the exception of
the Petersen graph, a connected vertex-transitive graph whose
automorphism group contains a transitive subgroup
with a cyclic commutator subgroup of prime-power order, is hamiltonian.
Furthermore, for connected vertex-transitive graphs of orders $p$, $2p$
(except for the Petersen graph), $3p$, $p^2$, $p^3$, $p^4$
and $2p^2$ it was shown that they are hamiltonian
(see \cite{BA79,CS94,YC98,DM85, DM87,DM88,T67}).
(Throughout this paper $p$ will always denote a prime number.)
On the other hand, connected vertex-transitive graphs
of orders $4p$ and $5p$
are only known to have Hamilton paths (see \cite{MP82,MP83}).
It is the object of this paper to complete the analysis of
hamiltonian properties of vertex-transitive graphs of order $4p$
by proving the following result.

\begin{theorem}
\label{the:main}
With the exception of the Coxeter graph,
every vertex-transitive graph of order $4p$, where $p$ is
a prime, is hamiltonian.
\end{theorem}

The proof of Theorem~\ref{the:main} is carried out over the
remaining sections.
Our strategy in the search for Hamilton
cycles in connected vertex-transitive graphs of order $4p$
is based on an analysis singling out
two facets of the structure of graphs in question.

First, a thorough analysis of various possibilities arising from
(im)primitivity of the action of the automorphism group of
a vertex-transitive graph of order $4p$ is done in Section~\ref{sec:pre1}.
More precisely, a vertex-transitive graph on $4p$ vertices
falls into one (but possibly more then one) of eight classes, depending
on various kinds of imprimitivity block systems its automorphism group admits
(see Table~\ref{tab:classes} in Section~\ref{sec:pre1} for details).
For some of these classes, sufficient conditions for existence
of Hamilton cycles in the corresponding graphs are given (see Lemmas \ref{lem:prim},
\ref{lem:classII} and \ref{lem:classIII}), leading to Proposition~\ref{pro:sec3},
where we prove that a connected vertex-transitive graphs of order $4p$
not isomoprhic to the Coxeter graph
is either hamiltonian or it has an imprimitivity block system with blocks of size $p$ or $2p$.

This result, reducing the total number of classes from
the initial eight to three, is then combined in Section~\ref{sec:four}
with results obtained from our second analysis taking into account
the well known fact that every vertex-transitive graph of order $mp$,
where $m\leq p$ has an $(m,p)$-semiregular automorphism \cite{DM81}.
In particular, letting  $\gamma$ be a $(4,p)$-semiregular
automorphism of a vertex-transitive graph $X$ of order $4p$,
the corresponding quotient graph $X_\gamma$ of $X$ with respect to $\gamma$
is one of six connected graphs of order $4$.
In \cite{MP83}, a thorough analysis for each of these six cases
resulted in the proof that every such graph has a Hamilton path.
Of course, as close as the concepts of Hamilton paths and cycles may seem,
the difficulties encountered in constructions of Hamilton cycles usually
greatly exceed those encountered in similar constructions of Hamilton paths.
It is therefore not surprising that this second approach alone
had not been enough to complete the result, thus calling for our two way analysis.


\section{Preliminary observations}
\label{sec:two}
\indent

Throughout this paper graphs are finite, simple, undirected
and connected, unless specified otherwise. By $p$ we shall always denote a prime number.
Also, all groups are assumed to be finite.
For adjacent vertices $u$ and $v$ in $X$,
we write $u \sim v$ and denote the corresponding edge by $uv$.
Given a graph $X$ we let $V(X)$, $E(X)$ and $\Aut X$
be the vertex set, edge set and the automorphism group of $X$, respectively.
A graph $X$ is said to be {\em vertex-transitive}
if its automorphism group $\Aut X$ acts transitively on $V(X)$.
Let $U$ and $W$ be disjoint subsets of $V(X)$.
The subgraph of $X$ induced by $U$ will be denoted by $X \la U \ra$;
in short, by $ \la U \ra$, when the graph $X$ is clear from the context.
Similarly, we let $X[U,W]$ (in short $[U,W]$) denote the bipartite subgraph
of $X$ induced by the edges having one endvertex in $U$
and the other endvertex in $W$.

Given a transitive group $G$ acting on a set $V$,
we say that a partition $\B$ of $V$ is $G$-{\em invariant}
if the elements of $G$ permute the parts, that is, {\em blocks} of $\B$, setwise.
If the trivial partitions $\{V\}$ and $\{\{v\}: v \in V\}$ are the only
$G$-invariant partitions of $V$, then $G$ is said to be {\em primitive},
and is said to be {\em imprimitive} otherwise.
In the latter case we shall refer to a corresponding
$G$-invariant partition as to an {\em imprimitivity block system} of $G$.
If the set $V$ above is the vertex set
of a vertex-transitive graph $X$,
and  $\B$ is an imprimitivity system of $G$, then clearly any
two blocks $B,B' \in \B$ induce isomorphic vertex-transitive
subgraphs.

For a graph $X$ and a partition $\PP$ of $V(X)$,
we let $X_{\PP}$ be the associated {\em quotient graph} of $X$ relative
to $\PP$, that is, the graph with vertex set $\PP$ and
edge set induced naturally by the edge set $E(X)$.
An automorphism of a graph is called
$(m,n)$-{\em semiregular}, where $m \geq 1$ and $n \geq 2$ are
integers, if it has $m$ orbits of length $n$ and no other orbit.
In the case when $\PP$ corresponds to
the set of orbits of a semiregular  automorphism $\gamma \in \Aut X$,
the symbol $X_{\PP}$ will be replaced by $X_{\gamma}$.

\begin{figure}[h!]
\begin{footnotesize}
\begin{center}
\includegraphics[width=0.40\hsize]{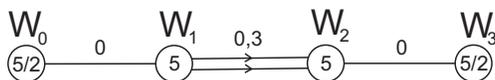}
\caption{\label{fig:example}\footnotesize The dodecahedron given in the Frucht's notation
relative to a $(4,5)$-semiregular automorphism.}
\end{center}
\end{footnotesize}
\end{figure}

Let $X$ be a connected vertex-transitive graph of order $4p$,
let $\W = \{W_i\ |\ i \in \ZZ_4\}$ be the set of orbits of a $(4,p)$-semiregular automorphism
$\gamma$ of $X$  and let the vertices of $X$ be labeled in such a way
that  $v_i^r\in W_i$ for $i\in\ZZ_4$ and $r\in \ZZ_p$.
Then $X$ may be represented by the notation of Frucht \cite{RF70}
emphasizing the four orbits of $\gamma$. (In fact Frucht's notation can be used for any graph that admits a
semiregular automorphism but we explain it here just for graphs admitting a $(4,p)$-semiregular
automorphism.) In particular, the  four orbits of $\gamma$ are  represented by four circles.
The symbol $p/x$, where $x\in\ZZ_p^*$, inside a circle coresponding to   the orbit $W_i$ means
that for each $r\in \ZZ_p$,
the vertex $v_i^r$ is adjacent to the vertex $v_i^{r+x}$. Similarly the symbol $p$
inside a circle coresponding to   the orbit $W_i$
means that $W_i$ is an independent set of vertices. Finally,
an arrow pointing from  the circle representing the orbit $W_i$ to the circle representing the
orbit $W_j$, $j\ne i$, labeled by $y\in\ZZ_p$ means that
for each $r\in\ZZ_p$, the vertex $v_i^r\in W_i$ is adjacent to the vertex $v_{j}^{r+y}$.
An example illustrating this notation is given in Figure~\ref{fig:example}.

The following classical result, due to Jackson \cite{J78}
giving a sufficient condition for the existence
of Hamilton cycles in regular graphs will be used here and throughout the rest of this paper.

 \begin{proposition}
 \label{pro:jack}
 {\rm \cite[Theorem~6]{J78}}
 Every $2$-connected regular graph of order $n$ and valency at least $n/3$
 is hamiltonian.
\end{proposition}

We end this section with the  proof that all vertex-transitive graphs of order $4p$,
$p \leq 5$ a prime, are hamiltonian. This will simplify the hamiltonicity analysis in the subsequent sections.
In the proof the so called LCF code \cite{RF77} will be used. The LCF code of a hamiltonian cubic
graph relative to one  of its Hamilton cycles $(v_0,v_1,\ldots, v_{n-1},v_0)$ is a list
$\LCF[a_0,a_1,\ldots,a_{n-1}]$ of elements of $\ZZ_n\setminus\{0,1,n-1\}$ such that $v_i$ is
adjacent to $v_{i+a_i}$ for every $i\in \ZZ_n$. In addition, if there exists a proper divisor $k$ of $n$
such that $a_i=a_{i+rk}$ for all $i\in\ZZ_k$ and $r\in\{1,2,\ldots,\frac{n}{k}-1\}$
then the notation is simplified to $\LCF[a_0,a_1,\ldots,a_{k-1}]^{\frac{n}{k}}$.

\begin{proposition}
\label{pro:less5}
A connected vertex-transitive graph of order $4p$, where $p\leq 5$ is a prime, is hamiltonian.
\end{proposition}

\begin{proof}
For $p=2$ the result follows from \cite{DM85}. By  \cite{McK79}
every vertex-transitive graph  of order $12$ is also a Cayley graph.
By Proposition~\ref{pro:jack}, it suffices to consider only graphs
of valency at most $3$. There are five such graphs: $C_{12}$,
$C_6\times K_2$, a graph obtained from $K_4$ by replacing each
vertex by a triangle, $\Cay(\ZZ_{12}, \{1,6\})$ and the graph with
LCF code $[5,-5]^6$. All of these graphs are hamiltonian. We may
therefore assume that $p=5$. (Note that by \cite{MR90}, there are
$1190$ connected vertex-transitive graphs of order $20$.) Let $X$ be
a connected vertex-transitive graph of order $20$. By
Proposition~\ref{pro:jack}, we may assume that the valency of $X$ is
less then $7$. Suppose first that  $X$ is a Cayley graph of a group
$G$ and let $P$ be a Sylow $5$-subgroup of $G$. Then $P$ is normal
in $G$ and the quotient group $G/P$, being of order $4$, is abelian.
Therefore, either $G$ itself is abelian or the commutator subgroup
of $G$ is cyclic of order $5$. Hence  by \cite{ED83,DM83} $X$ has a
Hamilton cycle. Let now $X$ be a non-Cayley graph of order $20$. It
can be deduced from \cite{MR90} that there are  $80$ possibilities
for $X$, with only $16$ having valency less than $7$. For these
graphs program package {\sc Magma} \cite{Mag} was used to find a
$(4,p)$-semiregular automorphism relative to which the
corrcorresponds to the dodecahedron which is known to possess a
Hamilton cycle. In all other cases  (with exception of the graph in
the second column of the third row, for which the existence of a
Hamilton cycle is straightforward) a Hamilton cycle is found using
the well known lifting of a Hamilton cycle in the quotient graph
(see also Proposition~\ref{pro:lem5}).

\end{proof}

\begin{figure}[h!]
\begin{footnotesize}
\begin{center}
\includegraphics[width=0.60\hsize]{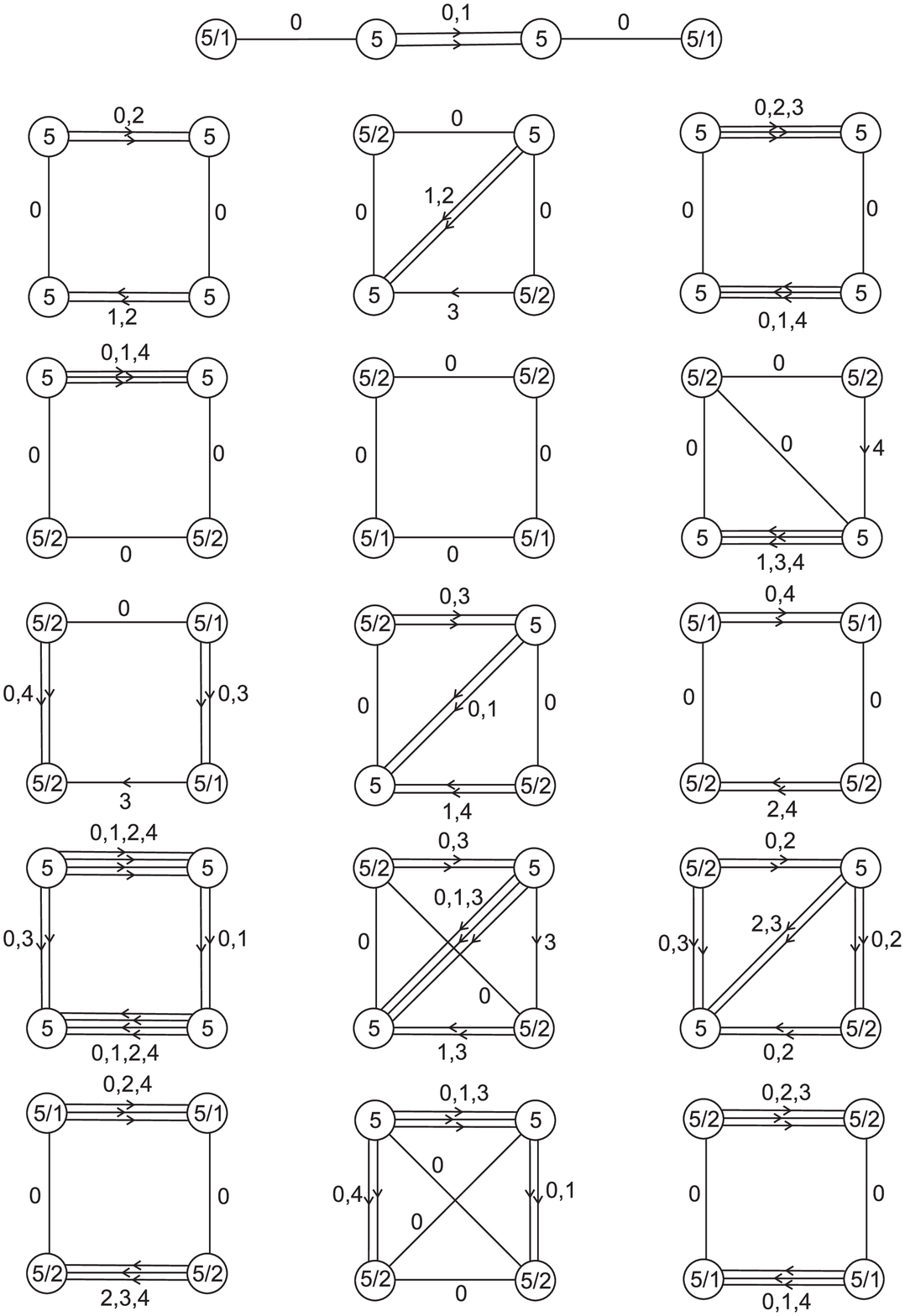}
\caption{\label{fig:graphs20}\footnotesize A list of all connected non Cayley vertex-transitive
graphs of order $20$ that are of valency less then $7$   given in the Frucht's notation.}
\end{center}
\end{footnotesize}
\end{figure}


\section{Analysis with respect to the action of $\Aut X$}
\label{sec:pre1}
\indent

\bigskip

An analysis of (im)primitivity of the full automorphism group
of a vertex-transitive graph of order $4p$, $p$ a prime,
is crucial in the proof of the main theorem of this paper.
Let us first divide all vertex-transitive graphs of order $4p$
into eight classes in the following way.
For a vertex-transitive graph $X$ of order $4p$,
let $A = \Aut X$ and choose $v \in V(X)$.
Let $(A_0,A_1,...,A_{k-1})$ be a sequence of groups such
that $A_0 =A$, $A_{k-1} =A_v$ is the vertex stabilizer and $A_i$
is maximal in $A_{i-1}$, $i \in \{1,...,k-1\}$. The corresponding
sequence of indices $[A_{i-1}:A_i]$, $(i\in \{1,...,k-1\})$, will be
called a {\em type} of the graph $X$. In view of these comments we
shall say that $X$ belongs to {\em Class~I}, {\em Class~II}, {\em
Class~III}, {\em Class~IV}, {\em Class~V}, {\em Class~VI}, {\em
Class~VII} and {\em Class~VIII}, respectively, if it is of type
$(4p)$, $(2:2p)$, $(2p:2)$, $(2:p:2)$, $(p:2:2)$, $(p:4)$, $(4:p)$
and $(2:2:p)$ (see also Table~\ref{tab:classes}).
For example, {\em Class~I} contains vertex-transitive graphs of order $4p$
with a primitive automorphism group and {\em Class~II} contains
vertex-transitive graphs of order $4p$ whose automorphism group has an imprimitivity
system of two blocks of size $2p$ and the block stabilizer acts
primitively on each of the two blocks. As we shall see in
Lemmas~\ref{lem:K=1} and \ref{lem:classV} the above eight
classes are not all disjoint.

\begin{table}[h!]
$$
\begin{array}{|l|c|}
\hline
{\it Class~I}&(4p)\\
\hline
{\it Class~II}&(2:2p)\\
\hline
{\it Class~III}&(2p:2)\\
\hline
{\it Class~IV}&(2:p:2)\\
\hline
{\it Class~V}&(p:2:2)\\
\hline
{\it Class~VI}&(p:4)\\
\hline
{\it Class~VII}&(4:p)\\
\hline
{\it Class~VIII}&(2:2:p)\\
\hline
\end{array}
$$
\caption{\label{tab:classes}{\footnotesize Eight classes of vertex-transitive graphs of order $4p$}}
\end{table}

\bigskip
The following result on primitive groups of degree $4p$ may be
extracted from \cite{LWWX, LS}.  By $D_{2n}$ we denote the dihedral group of order $2n$.

\begin{proposition}
\label{pro:3.1} Let $G$ be a primitive group of degree $4p$, where
$p\ge 7$, is a prime. Then $G$ is one of the
following:
\vspace{-5pt}
\begin{itemize}
\itemsep=0cm
\parskip=0pt\parsep=0pt
\item[(i)]  $A_8$ or $S_8$ acting on the $28=4p$ unordered pairs of points
    from an 8-element set;
\item[(ii)] $PSL(2,8)$ acting on the $28=4p$ cosets of a subgroup $D_{18}$;
\item[(iii)] $PGL(2,7)$ acting on the $28=4p$ cosets of a subgroup $D_{12}$;
\item[(iv)] $PSL(2,16)\le G\le P\Gamma L(2,16)$ acting on the $68=4p$
        cosets of a subgroup $N_G(PGL(2,4))$;
\item[(v)] $PSL(3,3)\le G\le PGL(3,3)$ acting on the $52=4p$ incident
    point-line pairs of $PG(2,3)$.
 \end{itemize}
\end{proposition}

Of course, vertex-transitive graphs  arising from  the above actions
in Proposition~\ref{pro:3.1} belong to {\it Class~I} and
{\sc Magma} program package \cite{Mag} was used
to obtain semiregular automorphisms relative to which
a Hamilton cycle in the corresponding quotient graph
lifting to a Hamilton cycle in the original graph was found.
It turns out that the Coxeter graph, a cubic graph associated
with the group action (iii), is the only graph not possessing a Hamilton cycle \cite{TRG}.
For details see Appendix~\ref{sec:appendix2}.
Combining the above arguments with Proposition~\ref{pro:less5} we  have the following result.

\begin{lemma}
\label{lem:prim}
Let $X$ be a connected vertex-transitive graph of order $4p$, $p$ a prime,
belonging to {\it Class~I}. Then $X$ is either hamiltonian or it is isomorphic to
the Coxeter graph.
\end{lemma}

The following result on primitive groups of degree $2p$ that may be deduced
from \cite{LS} will be needed here and later on in the paper.

\begin{proposition}
\label{pro:degree2p}
A primitive group $G$ of degree $2p$, $p$ a prime, is one of the following:\\
(i) $G$ is simply primitive and $p =5$
and $G = A_5$ or $G = S_5$;\\
(ii) $G = A_{2p}$ or $G = S_{2p}$;\\
(iii) $p =11$ and $G = M_{22}$;\\
(iv) $p = {\frac{1 + q^{2^t}}{2}}$, where $q$ is an odd prime,
$\Aut PSL(2,k)$ containing $PSL(2,k)$, where $k=q^{2^t}$ and $q$
is an odd prime. \\
Moreover, $G$ is simply primitive in case (i)
and is doubly transitive in all other cases.
\end{proposition}

For a permutation group $G$ acting on a set $V$ and a subset $W$
of $V$ we let $G_W$ denote the setwise stabilizer of $W$ in $G$ and we let
$G_{(W)}$ denote the pointwise stabilizer of $W$ in $G$.
The next two results assure the existence of Hamilton cycles in vertex-transitive graphs
of order $4p$ belonging to {\it Classes~II} and {\it III}.

\begin{lemma}
\label{lem:classII}
A connected vertex-transitive graph of order $4p$, $p$ a prime,
belonging to {\it Class~II} is hamiltonian.
\end{lemma}

\begin{proof}
By Proposition~\ref{pro:less5}, we may assume that $p\ge 7$.
Let $X$ be  a connected vertex-transitive graph with $4p$ vertices,
let  $A=\Aut X$ be its automorphism group, and let
$\B=\{B,B'\}$ be an imprimitivity block system of $A$ consisting
of two blocks of size $2p$. Since $X$ is of type $(2:2p)$, the
group $A_{B}=A_{B'}$ is a primitive group of degree $2p$, in its action
on $B$ and $B'$. Now, in view of Proposition~\ref{pro:degree2p},
these two actions are equivalent and $A_{B}=A_{B'}$ acts doubly transitively on $B$ and $B'$.
For regularity reasons, the induced subgraphs
on $B$ and $B'$ are either both isomorphic to the complete graph
$K_{2p}$ or are totally disconnected. In the first case,
the valency of $X$ is greater then
$2p-1$, and hence $X$ is hamiltonian by Proposition~\ref{pro:jack}.
If $X\la B \ra$ and $X\la B' \ra$ are totally disconnected
then, depending on whether the two actions are faithful or unfaithful,
we obtain that  either $X\cong K_{2p,2p}-2pK_2$ or $X \cong K_{2p,2p}$.
In both cases, Proposition~\ref{pro:jack} gives us a Hamilton cycle in $X$.

\end{proof}

\begin{lemma}
\label{lem:classIII}
A connected vertex-transitive graph of order $4p$,  $p$ a prime,
belonging to {\it Class~III} is hamiltonian.
\end{lemma}

\begin{proof}
By Proposition~\ref{pro:less5}, we may assume that $p\ge 7$.
Let $X$ be a connected vertex-transitive graph with $4p$ vertices,
let  $A=\Aut X$ be its automorphism group, let
$\B$ be an imprimitivity block system of $A$ consisting
of $2p$ blocks of size $2$, and let $K$ be the kernel
of the action of $A$ on $\B$.
Since $X$ is of type $(2p:2)$,
it follows that $\bar{A} =A/K$ is primitive on $\B$.
By Proposition~\ref{pro:degree2p},  $\bar{A}$ acts doubly transitively and so
the quotient graph $X_\B$ is isomorphic to the complete graph $K_{2p}$.
Now, the bipartite subgraphs $X[B,B']$, $B, B' \in \B$,
are all isomorphic and must have, for arithmetic reasons, an even number of edges.
Hence $X[B,B'] \cong 2K_2$ or $X[B,B'] \cong K_{2,2}$. Therefore
the valency of the graph $X$ is at least $2p-1$
and Proposition~\ref{pro:jack} gives us a Hamilton cycle in $X$.

\end{proof}

\begin{lemma}
\label{lem:K=1}
Let $X$ be a connected vertex-transitive graph of order $4p$, $p\geq 7$ a prime,
such that  $\Aut X$ admits an imprimitivity block system $\B$ with $p$
blocks of size $4$ (and so $X$ is either in {\it Class~V} or {\it Class~VI}).
If the kernel $K$ of the action of
$\Aut X$ on $\B$ is trivial then $X$ belongs to
{\it Class~IV}, {\it Class~VII} or {\it Class~VIII}, or $X$ is isomorphic to the graph shown in
Figure~\ref{fig:PSL(3,2)}.
\end{lemma}

\begin{proof}
Since, by assumption, $K=1$ we have that $A=\Aut X \cong \bar{A}=A/K$ is a group of prime degree.
If $A$ is solvable then, in view of \cite[Proposition~2.1]{MSZ95}, we have that $A\le A(1,p)$ and it
follows from \cite[Theorem~3.5B]{DM96} that $A$ has a regular normal Sylow $p$-subgroup.
Thus, there exists a $(4,p)$-semiregular element $\gamma\in A$ such that
$\langle \gamma \rangle$ is normal in $A$ and so, by \cite[Theorem~8.8]{W64}, $X$ belongs
to {\it Class~VII} or {\it Class~VIII}.

Suppose now that $A$ is nonsolvable. Then, by \cite[Theorem~3.5B]{DM96},
$A$ is doubly transitive and so  $X_\B=K_p$.  Again
using Proposition~\ref{pro:degree2p} and checking all the possibilities for the
existence of index 4 subgroups in the block stabilizer $A_{B}$, $B\in {\cal B}$,
we can see that $PSL(n,k)\le A\le \Aut PSL(n,k)$ for appropriate $n$ and $k$,
in view of the fact that  $p \geq 7$.

If $A=PSL(n,k)$ or if $A$ properly contains a copy of $PSL(n,k)$ acting transitively,
then following the argument used in \cite{MSZ95} we obtain
that the groups $PSL(3,2)$ and $PSL(3,3)$ acting on cosets of $S_3$ and $2S_3$,
respectively, are the only possibilities. The latter is clearly impossible for
it would give rise to a graph of order $468 = 4 \cdot 117$, which is not of the form $4p$.
As for the action of $PSL(3,2)$ on $S_3$, using program package {\sc Magma} \cite{Mag} we deduce that
$S_3$ has six nontrivial suborbits, two of which are non-self-paired of length $6$.
Of the four self-paired suborbits, three are of length $3$ and one is of length $6$.
The graph arising from the union of the two non-self-paired suborbits
has valency $12$ and is isomorphic to the graph arising from the self-paired suborbit
of length $12$ in the action of  $PGL(2,7)$ on cosets of $D_{12}$.
The graph arising from one of the  suborbits
of length $3$ is isomorphic to the Coxeter graph and hence with a primitive automorphism group.
Next, the graphs arising from the other two suborbits of length $3$ are both disconnected and isomorphic to $7K_4$.
Furthermore, the union of these two graphs is isomorphic to the graph arising from one
of the self-paired suborbits of length $6$ in the action of  $PGL(2,7)$ on the cosets of  $D_{12}$.
As for the  graph arising from the union of two self-paired suborbits of length $3$,
one giving  rise to $7K_4$ and the other giving to the Coxeter graph, it is isomorphic
to the graph depicted in Figure~\ref{fig:PSL(3,2)} using Frucht's notation \cite{RF70}.
Finally, the graph arising from the self-paired suborbit  of length $6$ is isomorphic to one of the graphs
associated with the action of $PGL(2,7)$ on cosets of $D_{12}$.

If $A$ properly contains a copy of $PSL(n,k)$ acting intransitively, then
the normality of $PSL(n,k)$ in $A$ gives us an imprimitivity block system $\cal{C}$ for $A$.
Since $p$ does not divide $[\Aut PSL(n,k):PSL(n,k)]$,
it follows that $\cal{C}$ consists  of blocks of length $p$ or $2p$, completing the proof of Lemma~\ref{lem:K=1}.

\end{proof}

\begin{figure}[ht!]
\begin{footnotesize}
\begin{center}
\includegraphics[width=0.20\hsize]{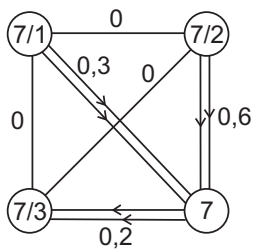}
\caption{\label{fig:PSL(3,2)}\footnotesize The vertex-transitive graph whose automorphism group
is isomorphic to $PSL(3,2)$ given in the Frucht's notation relative to a $(4,7)$-semiregular automorphism.}
\end{center}
\end{footnotesize}
\end{figure}

Vertex-transitive graphs of order $2p$, $p$ a prime,
were described  in \cite{DM81}.
Among others it was proved there that, provided a vertex-transitive graph $X$
of order $2p$ admits an imprimitve group $G$
(with blocks of size $p$ or $2$), one can always find an imprimitive
subgroup of $G$ which has blocks of size $p$.
Moreover, if $A = \Aut X$ itself has blocks of size $2$ and no blocks of size $p$,
it may be deduced from the proof of \cite[Theorem 6.2]{DM81} that
$X$ or its complement is the wreath product $Y \wr 2K_1$ where $Y$ is a  $p$-circulant
(Recall that for graphs $X$ and $Y$, the {\em wreath product}, sometimes also called
the {\em lexicographic product} $X\wr Y$,
has vertex set $V(X)\times V(Y)$ with two vertices $(a,u)$ and $(b,v)$
adjacent in  $X\wr Y$  if and only if either  $ab\in E(X)$ or  $a=b$ and $uv\in E(Y)$.).
This enables us to prove the following result.

\begin{lemma}
\label{lem:classV}
Let $X$ be a connected vertex-transitive graph of order $4p$, $p\geq 7$,
belonging to {\it Class~V} or {\it Class~VI} and let $\B$ be an imprimitivity block system of
$\Aut X$ with blocks of size $4$. Then one of the following holds
\vspace{-5pt}
\begin{enumerate}[(i)]
\itemsep=0cm
\parskip=0pt\parsep=0pt
\item  $X$  belongs to {\it Class~IV}, {\it Class~VII} or {\it Class~VIII};
\item $X$ is as shown in Figure~\ref{fig:PSL(3,2)};
\item  $X$  is a Cayley graph of an abelian group;
\item $X$ is isomorphic  to $Y\wr Z$, where $Y$ is a connected vertex-transitive graph of order $2p$
and $Z$ is either $2K_1$ or $K_2$;
\item $X$ is a regular $\ZZ_2$-cover of $K_p\wr 2K_1$; or
\item  there exist adjacent blocks $B,B'$ in $X_\B$  such that $X[B,B']$ is $K_{4,4}$ or $2C_4$.
\end{enumerate}
\end{lemma}

\begin{proof}
Let $K$ be the kernel of the action of $A=\Aut X$ on $\B$.
If $K=1$, then Lemma~\ref{lem:K=1} implies that  $X$  belongs to {\it Class~IV}, {\it Class~VII} or {\it Class~VIII},
or $X$ is the graph in Figure~\ref{fig:PSL(3,2)}.
Assume now that $K$ is nontrivial.  We shall distinguish two different cases.

\medskip

\noindent
{\sc Case~1.}
If $K$ is intransitive on each of the blocks in ${\cal B}$, it follows that
$K^B$ is either $\ZZ_2$ for each $B\in{\cal B}$ or $\ZZ_2^2$ for each $B\in{\cal B}$,
and further, the orbits of $K$ form an imprimitivity block system ${\cal E}$ with blocks of size $2$.
Clearly, $K$ is also the kernel of the action of $A$ on ${\cal E}$.
If $K \neq \ZZ_2$ then the action of $K$ on the  blocks in ${\cal E}$
is unfaithful and so $X$ must be the  wreath product
of the vertex-transitive graph $X_{\cal E}$ of order $2p$ with $2K_1$ or with $K_2$,
and so (iv) holds.
So let $K = \ZZ_2$. Consider the group $\bar{A} = A/K$
acting on ${\cal B}$. If $\bar{A}$ is solvable, then it has a normal
subgroup $PK/K$ of order $p$ where $P$ is a Sylow $p$-subgroup of $A$.
Since $K = \ZZ_2$ the Sylow theorems imply that $P$ is a characteristic subgroup
of $PK$. Since $PK$ is normal in $A$ we have that $P$ is normal in $A$.
It follows that $X$ belongs to {\it Class VII} or {\it Class VIII}.
We may therefore  assume that the action of $\bar{A}$ on ${\cal B}$ is unsolvable and hence doubly transitive, by
Burnside's classical result (see \cite[Theorem~7.3]{Pa68}).
Hence $X_\B = K_p$. Consider the action of $\bar{A}$ on the quotient graph
$X_{\cal E}$. If apart from blocks of size $2$
it also has blocks of size $p$, then $X$ belongs to {\it Class IV}.
So we may assume that $\bar{A}$ as well as $\Aut X_{\cal E}$ has no blocks of size $p$.
By the comments preceding the statement of Lemma~\ref{lem:classV} and taking into account
the fact that $X_{\cal B} = K_p$, it follows that  $X_{\cal E}$
is isomorphic to the wreath product $K_p \wr 2K_1$. Consequently,  $X$ is isomorphic either
to $X_{\cal E} \wr 2K_1$ or to $X_{\cal E} \wr K_2$, or it is a regular $\ZZ_2$-cover of $X_{\cal E}$.
In short,  either (iv) or (v) holds.

\medskip

\noindent
{\sc Case~2.} Assume now that $K$ is transitive on each of the blocks  $B\in\B$.
We  have $K^B \in\{\ZZ_2^2,\ZZ_4,D_8,A_4,S_4\}$.
Suppose first that $K$ is faithful. Then $K\in\{\ZZ_2^2,\ZZ_4,D_8,A_4,S_4\}$
and we can assume that there is a characteristic subgroup $H$ in $K$ of order 4
(either $\ZZ_2^2$ or $\ZZ_4$). Hence $H$ is normal in $A$ and so $H$ is normal in $\langle \gamma,H\rangle $,
where $\gamma$ is some $(4,p)$-semiregular element in $A$.
The Sylow theorems imply  that $\langle \gamma,H\rangle =H\times \langle \gamma\rangle $.
Hence $X$ is a Cayley graph either of $\ZZ_{4p}$ or of $\ZZ_{2p}\times \ZZ_2$,
and so (iii) holds.

We may now assume that $K$ is unfaithful. Let $B, B' \in {\cal B}$ be adjacent
in $X_\B$. Then $K_{(B)}^{B'}\ne 1$ and $K_{(B)}^{B'}$ is normal in $K^{B'}$.
If $K_{(B)}^{B'}$ is transitive then $X[B,B']=K_{4,4}$.
If $K_{(B)}^{B'}$ is intransitive then  clearly $K^{B'} = K^B \in \{\ZZ_2^2, \ZZ_4,D_8\}$.
Moreover,  $K_{(B)}^{B'}$  must have two orbits on $B'$ and either $X[B,B'] =K_{4,4}$
or $X[B,B']=2C_4$, and so (vi) holds. This completes the proof of Lemma~\ref{lem:classV}.

\end{proof}

Given a graph $X$ admitting a $(4,p)$-semiregular automorphism with the set of orbits $\W$
and an imprimitivity block system  $\B$ of $\Aut X$, we have that

\begin{equation}
\label{orbite-bloki}
|W \cap B|=1 \; \; {\rm or}\; \; W\subseteq B,
\end{equation}

\noindent
for each $W\in \W$ and $B\in \B$.

Combining together results of this section we can prove the following proposition that reduces
the possible existence of nonhamiltonian graphs of order $4p$ to {\it Class~IV},
{\it Class~VII} or {\it Class~VIII}.

\begin{proposition}
\label{pro:sec3}
Let $X$ be a connected vertex-transitive graph of order $4p$, $p$ a prime
not isomoprhic to the Coxeter graph. Then either $X$ is
hamiltonian or $X$ belongs to  {\it Class~IV}, {\it Class~VII} or {\it Class~VIII}. In short,
either $X$ is hamiltonian or $\Aut X$ has an imprimitivity block system with blocks of size $p$ or $2p$.
\end{proposition}

\begin{figure}[ht!]
\begin{footnotesize}
\begin{center}
\includegraphics[width=0.65\hsize]{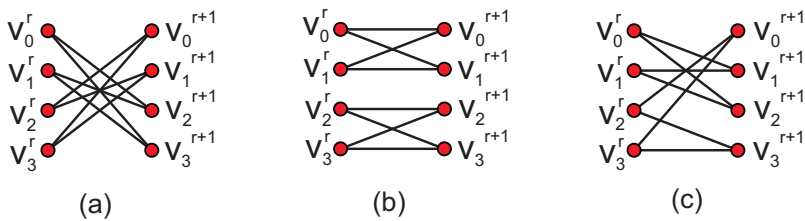}
\caption{\label{fig:2C4}\footnotesize Possible forms of the bipartite graph $X[B,B']$ where $B$ and $B'$
are adjacent blocks of size $4$.}
\end{center}
\end{footnotesize}
\end{figure}

\begin{figure}[ht!]
\begin{footnotesize}
\begin{center}
\includegraphics[width=0.45\hsize]{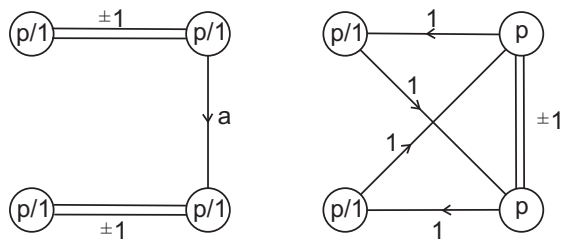}
\caption{\label{fig:fruc}\footnotesize Two possibilities for a spanning subgraph in $X$. The graph on the left
coresponds to the graph in Figure~\ref{fig:2C4}(b) and the graph on the right
coresponds to the graph in Figure~\ref{fig:2C4}(c). Where $a\in \ZZ_p$.}
\end{center}
\end{footnotesize}
\end{figure}

\begin{proof}
By Proposition~\ref{pro:less5}, we may assume that $p\ge 7$.
Let $X$ be a connected vertex-transitive graph of order $4p$ that
belongs  to {\it Class~I}, {\it Class~II}, {\it Class~III}, {\it Class~V} or {\it Class~VI}.
Then  by Lemmas~\ref{lem:prim}, \ref{lem:classII}  and \ref{lem:classIII},
we may assume that $X$ belongs to {\it Class~V} or {\it Class~VI}.
Then  one of the statements (ii)-(vi) in Lemma~\ref{lem:classV} holds.

First, if (ii) holds and  $X$ is as in Figure~\ref{fig:PSL(3,2)}, then clearly $X$ is hamiltonian.
Also, if (iii) holds then $X$ is hamiltonian too, in view of \cite{DM83}.
If (iv) holds, then $X$ is hamiltonian in view of the fact that
a connected vertex-transitive graph of order $2p$, $p\geq7$, has a Hamilton cycle  \cite{BA79}),
and the fact that the   wreath product of a hamiltonian graph with  $2K_1$ is hamiltonian.
If (v) holds and $X$ is a regular $\ZZ_2$-cover of $K_p\wr 2K_1$,
then its valency is $2p-2$, and hence,  by Proposition~\ref{pro:jack}, $X$ is hamiltonian.
Finally, let us assume that (vi) holds.
If there  there exist adjacent blocks $B,B'$ in $X_\B$ such that $X[B,B']$ is isomorphic to $K_{4,4}$
then $X$ is clearly hamiltonian.
We may therefore assume that for any two adjacent blocks
 $B,B'$ in $X_\B$ the graph $X[B,B']$ is isomorphic to $2C_4$.

 Let $\W = \{W_i\ |\ i \in \ZZ_4\}$ be the set of orbits of a $(4,p)$-semiregular automorphism
$\gamma$ of $X$. By (\ref{orbite-bloki}), there are $v_0\in W_0$, $v_1\in W_1$,  $v_2\in W_2$ and
$v_3\in W_3$ such that $B=\{v_0,v_1,v_2,v_3\}$ is a block.  Let $v_i^r=\gamma^r(v_i)$, for $i\in\ZZ_4$
and $r\in\ZZ_p$. Then $\B=\{B_r\mid r\in \ZZ_p\}$ where $B_r = \gamma^r(B)=\{v_o^r,v_1^r,v_2^r,v_3^r\}$
for $r\in\ZZ_p$.
Without loss of generality we may assume that the bipartite graph $X[B,B']$ is
one of the graphs in Figure~\ref{fig:2C4}.
Now, the bipartite graph $X[B,B']$ in Figure~\ref{fig:2C4}(a)
gives rise to a spanning subgraph in $X$ that is isomorphic to
the wreath product of a connected vertex-transitive graph of order $2p$ with $2K_1$.
Clearly, in this case $X$ is  hamiltonian.
We may therefore assume that  $X[B,B']$
is either the one in  Figure~\ref{fig:2C4}(b) or the one in Figure~\ref{fig:2C4}(c).
It follows that  $X$ contains a  spanning subgraph isomorphic, respectively, to the graphs shown in  Figure~\ref{fig:fruc},
using Frucht's notation \cite{RF70}, with  $a \in \ZZ_p$.
Since

$$
v_1^0v_1^{-1}\cdots v_1^lv_3^{l+1}v_3^{l+2}\cdots v_3^{l}v_1^{l-1}v_1^{l-2}\cdots
v_1^2v_1^1v_2^{a+1}v_2^{a+2}\cdots
$$
$$
\cdots v_2^{a+k}v_0^{a+k+1}v_0^{a+k+2}\cdots v_0^{a+k}v_2^{a+k+1}v_2^{a+k+2}\cdots v_2^av_1^0
$$
is a Hamilton cycle in the graph on the left in Figure~\ref{fig:fruc} and
$$
v_0^0,v_2^1v_3^2v_1^4v_0^5v_2^6\cdots v_0^{p-4}v_2^{p-3}v_3^{p-2}v_1^{p-1}v_0^0
$$
is a Hamilton cycle in the graph on the right  in Figure~\ref{fig:fruc}, the result follows.

\end{proof}

\section{Analysis with respect to the quotient graph $X_{\gamma}$}
\label{sec:four}
\indent

\bigskip

\begin{figure}[ht!]
\begin{footnotesize}
\begin{center}
\includegraphics[width=0.80\hsize]{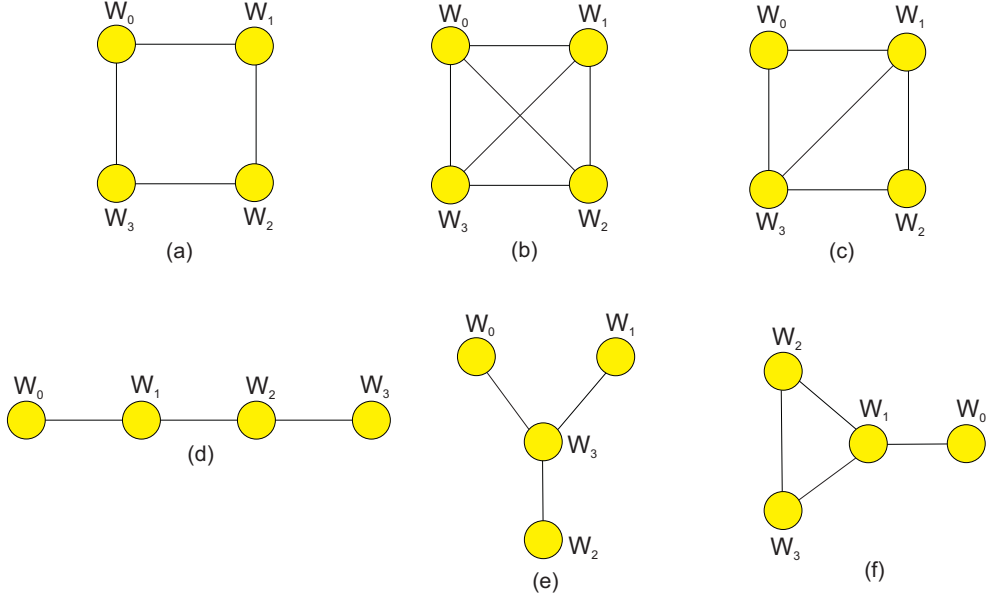}
\caption{\label{forms}\footnotesize The six possibilities for the quotient
graph $X_{\gamma}$ of a connected vertex-transitive graph $X$ of order $4p$.}
\end{center}
\end{footnotesize}
\end{figure}

We shall now combine Proposition~\ref{pro:sec3} with an analysis
of the quotient graph $X_{\gamma}$ of a connected vertex-transitive graph
$X$ of order $4p$, $p\geq7$, relative to a $(4,p)$-semiregular automorphism $\gamma$
which exists in $X$ in view of \cite[Theorem 3.4.]{DM81}.
Let $\W=\{W_i \mid i\in\ZZ_4\}$ denote the set of orbits of $\gamma$.
Now there are six different possibilities for the
quotient graph $X_{\gamma}$  of $X$ relative to $\gamma$ (see Figure~\ref{forms}).

The following easy observations are straightforward.
First, for any orbit $W_i$ of $\gamma$, the induced subgraph
$\la W_i \ra$ is regular of some even valency $d(W_i)$.
Moreover, if $d(W_i) > 0$, then $\la W_i \ra$ contains
a Hamilton cycle. Second,
for distinct $i,j$, the bipartite graph
$X[W_i,W_j]$ is regular of some valency
$d[W_i,W_j]\ge 0$. And finally, when
$d[W_i,W_j]\ge 2$, $X[W_i,W_j]$ contains a Hamilton cycle.

A graph is {\it Hamilton-connected} if for every pair of vertices
$u$ and $v$ there exists a Hamilton path whose endvertices are
$u$ and $v$. The following three results taken, respectively from,
\cite[Theorem~4]{CQ81}, \cite[Lemma~5]{MP82}, \cite[Theorem~3.9]{BA89},
will play an essential role in the proof of Theorem~\ref{the:main}.

\begin{proposition}
\label{pro:alspach}
(\cite{CQ81}). Let $X$ be a connected Cayley graph of
an abelian graph of valency at least $3$. If $X$ is not bipartite
then $X$ is Hamilton-connected.
\end{proposition}

\begin{proposition}
\label{pro:lem5}
(\cite{MP82}).
Let $\gamma$ be a semiregular automorphism of a graph $X$ and let
$C=W_0W_1\cdots W_{k-1}$, $k\geq3$, be a cycle in $X_{\gamma}$.
If $\la C \ra$ does not contain a Hamilton cycle,
then $d[W_i,W_{i+1}] =1$ for $i\in \ZZ_k$, and the graph induced by
the edges of the graphs $[W_i,W_{i+1}]$, $i\in \ZZ_k$, is a
disjoint union of $p$ cycles of length $k$ in $X$.
\end{proposition}

For the third result we need the concept of a coil of a cycle in a
quotient graph, introduced in \cite{BA89}. Let $X$ be a graph that
admits an $(m,n)$-semiregular  automorphism  $\alpha$ and let
$\W=\{W_i\mid i\in\ZZ_m\}$ be the set of orbits of $\alpha$. Let
$C=W_rW_sW_t\ldots W_qW_r$ be a cycle of length $k$ in $X_\W$ and
let $v_r^0$, $v_r^1$, $\ldots$, $v_r^{n-1}$ be a cyclic labelling of
the vertices of $W_r$ under the action of $\alpha$. Consider the
path of $X$ arising from a lifting of $C$, namely, start at $v_r^0$
and choose an edge from $v_r^0$ to a vertex $v_s^a$ of $W_s$. Then
take an edge from $v_s^a$ to a vertex of the $W_t$ following $W_s$
in $C$. Continue this way until returning to a vertex $v_r^b$ of
$W_r$. If $b\ne 0$, a path of length $k$ has been constructed and if
$b=0$, it is a cycle of length $k$. There will be more then one such
path if the degree between two consecutive orbits of $\alpha$ is
larger then one. The set of all paths in $X$ arising from a lifting
of $C$ is denoted by $coil(C)$. The following result is proved in
\cite{BA89}.

\begin{proposition}
\label{the:alspach}
(\cite{BA89}). Let $X$ be a graph admitting
an $(m,n)$-semiregular automorphism $\alpha$,
with $m \geq 4$ even and $n\geq3$,
and let $\W=\{W_i \mid i\in\ZZ_m\}$ be the set of orbits  of $\alpha$
such that each $\la W_i \ra$ has valency $2$ and is connected.
If $X_\W$ contains a Hamilton cycle $C$ such that $coil(C)$ contains a cycle,
then $X$ has a Hamilton cycle.
\end{proposition}

We are now ready to prove Theorem~\ref{the:main}.

\bigskip

\noindent \proofT \ Let $X$ be a connected vertex-transitive graph
of order $4p$ and valency $d=d(X)$, different from the Coxeter
graph. By Proposition~\ref{pro:less5}, we may assume that $p\ge 7$.
Moreover, we may also assume that $d\geq3$. Let $\W = \{W_i\ |\ i
\in \ZZ_4\}$ be the set of orbits of a $(4,p)$-semiregular
automorphism $\gamma$ of $X$. For $i \in \ZZ_4$, let $d_i$ denote
the valency of the induced graph $\la W_i \ra$, and for $i,j \in
\ZZ_4$, let $d_{i,j}$ denote the valency of the induced bipartite
graph $[W_i,W_j]$. By Proposition~\ref{pro:sec3}, we may assume that
$A = \Aut X$ has an imprimitivity block system $\B$ with blocks of
size $p$ or $2p$.


\begin{case} \end{case}
$X_{\gamma}$ has a $4$-cycle $W_0W_1W_2W_3W_0$
(see  Figure~\ref{forms}{{\small (a),(b),(c)}}).

\smallskip

\noindent
By Proposition~\ref{pro:lem5}, we may assume that $d_{i,i+1}=1$ for $i \in \ZZ_4$
and that the subgraph of $X$ spanned by all the edges of the graphs $[W_i,W_{i+1}]$, $i\in\ZZ_4$,
is a disjoint union of $p$ cycles of length $4$.


\begin{subcase}\end{subcase}  $X_{\gamma}$ is  the $4$-cycle $C_4$.

\smallskip

\noindent
The connectedness and regularity of
$X$ imply that $d_i=d-2 \ge 2$ for $i \in \ZZ_4$.
If $d_i=2$ for $i\in\ZZ_4$, then by Proposition~\ref{the:alspach}, $X$ has a Hamilton cycle.
If on the other hand, $d_i \ge 4$ for $i\in \ZZ_4$, Proposition~\ref{pro:alspach}
implies that each subgraph $\la W_i \ra$ is Hamilton-connected, and consequently $X$
is hamiltonian.


\begin{subcase}\end{subcase} $X_{\gamma}$ is the complete graph $K_4$.

\smallskip

\noindent By Proposition~\ref{pro:lem5}, we may assume that
$d_{i,j}=1$ for distinct $i,j \in \ZZ_4$. By regularity of $X$ we
have that $d_i =d_j$ for $i,j \in \ZZ_4$. If $d_i = 2$ for $i \in
\ZZ_4$, Proposition~\ref{the:alspach} implies that $X$ contains a
Hamilton cycle. On the other hand, if $d_i \ge 4$ for $i \in \ZZ_4$,
then Proposition~\ref{pro:alspach} implies that $X$ is hamiltonian.


\begin{subcase}\end{subcase} $X_{\gamma}$ is neither $C_4$ nor $K_4$
(see Figure~\ref{forms}{{\small(c)}}).

\smallskip

\noindent We may assume that $d_{1,3}>0$ and $d_{0,2}=0$. Therefore
the valency of vertices in $\la W_0 \cup W_2\ra$ is even, and so
$d\ge 4$ is even. Hence $d_0 = d_2\ge 2$. Consequently, $d_{1,3}$ is
even too, and so $d_{1,3} \geq 2$. Using (\ref{orbite-bloki}), it is
easily seen that $\B$ cannot consist of four blocks of size $p$ and
so it consists of two blocks of size $2p$, each being a union of two
orbits of $\gamma$. Without loss of generality we may assume that
either $W_0\cup W_1$ or $W_0\cup W_2$ is a block $B$ in $\B$. But
the former cannot occur as then $\la B \ra$ is not a regular graph.
If however the latter is the case, then $B'= W_1\cup W_3$ is the
other block in $\B$ inducing a connected graph, whereas $\la B \ra$
is disconnected, a contradiction.


\begin{case}
\label{drevo}
\end{case}
$X_{\gamma}$ is a tree.


\begin{subcase}
\end{subcase}
$X_{\gamma}$ is the $3$-path (see Figure~\ref{forms}{{\small(d)}}).

\smallskip

\noindent  By regularity, $d_0,d_3, d_{1,2} \ge 2$.
Assume first that $\B$ is an imprimitivity block system consisting of four blocks of
size $p$. By (\ref{orbite-bloki}),
$\B$ coincides with the set of orbits $\W$ of $\gamma$.
Since any two blocks give rise to
isomorphic vertex-transitive graphs,
it follows that $d_i =d _j$ for $i,j\in \ZZ_4$.
But then, as $d_{1,2} \ge 2$, the vertices in $W_1 \cup W_2$
would be of greater valency than those in $W_0 \cup W_3$, a contradition.

Assume now that $\B$ is an imprimitivity block system with two
blocks of size $2p$. By (\ref{orbite-bloki}) each block in $\B$ is a
union of two orbits of $\gamma$. In particular one of the sets  $W_0
\cup W_1$ or $W_0 \cup W_2$ or $W_0 \cup W_3$ must be a  block in
$\B$. The first possibility cannot occur for obvious arithmetic
raesons since $d_0 \neq d_1$. The second possibility implies that
$d_0 =d_2$ and $d_1 =d_3$. But then, on the one hand, comparing the
valencies of vertices in $W_0$ and $W_1$, it follows that $d_0 - d_1
= d_{1,2} \ge 2$, and on the other hand, comparing the valencies of
vertices in $W_2$ and $W_3$, it follows that $d_1 - d_0 = d_{1,2}
\ge 2$, a contradiction. Finally, the third possibility is also
impossible as $W_0 \cup W_3$ induces a disconnected graph, but $W_1
\cup W_2$ induces a connected graph.


\begin{subcase}
\end{subcase}
$X_{\gamma}$ is the star $K_{1,3}$ (see Figure~\ref{forms}{{\small(e)}}).

\smallskip

\noindent  By regularity, $d_3$ is clearly different (smaller) from
each of $d_i$, $i\in \ZZ_4\setminus \{3\}$. In particular, in view of (\ref{orbite-bloki}),
this implies that $\B$ does not consist of four blocks of size $p$.
Hence, $\B$ consists of two blocks of size $2p$.
By (\ref{orbite-bloki})
each block in $\B$ is a union of two orbits of $\gamma$.
Without loss of generality these two blocks are $W_0 \cup W_1$ and $W_2 \cup W_3$.
But the latter induces a graph which is not regular, a contradiction.


\begin{case}
\end{case}
$X_{\gamma}$ is the graph shown in
Figure~\ref{forms}{{\small (f)}}.

\smallskip

\noindent By regularity, $d_0 \ge 2$ and $d_1 \neq d_0$. This
implies that $A$ cannot have blocks of size $p$ , and so, using
(\ref{orbite-bloki}) again, $\B$ consists of two blocks of size
$2p$, each a union of two orbits of $\gamma$. For regularity reasons
$W_0 \cup W_1$ cannot be a block, and so with no loss of generality
the blocks must be $W_0 \cup W_2$ and $W_1 \cup W_3$. In particular
$d_0 = d_2$ and $d_1 = d_3$. But then, on the one hand, comparing
the valencies of vertices in $W_0$ and $W_2$, it follows that
$d_{0,1} = d_{1,2} + d_{2,3}$, and on the other hand, comparing the
valencies of vertices in $W_1$ and $W_3$, it follows that $d_{2,3} = d_{1,2} + d_{0,1}$,
a clear contradiction. This completes the proof
of Theorem~\ref{the:main}. 
\begin{flushright}\Qed\end{flushright}


\section{Appendix -- Vertex-transitive graphs from  Class~I}
\label{sec:appendix2}
\indent

 We discuss here hamiltonicity properties
 of vertex-transitive graphs of order $4p$ and valency less than $4p/3$
  having a primitive  automorphism group,  and thus
 arising from the actions in Proposition~\ref{pro:3.1}. 
 The graphs are given in  Table~\ref{tab:simbol} using a certain collection of subsets
 of $\ZZ_p$ associated with a $(4,p)$-semiregular automorphism.
 
\begin{figure}[ht!]
 \begin{footnotesize}
 \begin{center}
 \includegraphics[width=0.24\hsize]{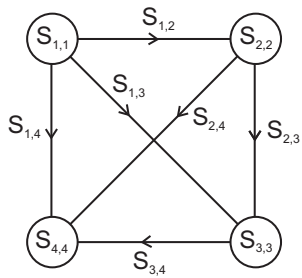}
 \caption{\label{fig:symbol-frucht}\footnotesize The Frucht's notation of a graph with symbol
 ${\bf S} =(S_{i.j})$.}
 \end{center}
 \end{footnotesize}
 \end{figure}

 Given a graph $X$ with a $(4,p)$-semiregular automorphism $\gamma$
 with orbits $W_i$, $i \in \ZZ_4$,
 choose $w_i \in W_i$ and define the following subsets of $\ZZ_p$,
 the collection of which determines $X$ uniquely.
For $i,j \in \ZZ_4$, we let 
 $S_{i,j} =\{s \in \ZZ_p: [w_i,\gamma^s w_j] \in E(X)\}$.
 Clearly $S_{j,i} =-S_{i,j}$.
 The $4\times 4$-"matrix" ${\bf S} =(S_{i,j})$
 whose $(i,j)$-th entry is
 the set $S_{i,j}$ is usually referred to as the {\em symbol} of $X$ relative to $\gamma$.
 The connection between the symbol of a graph that admits a $(4,p)$-semiregular automorphism
  and the Frucht's notation \cite{RF70} of a graph is given in Figure~\ref{fig:symbol-frucht}.
 
 As remarked in Section~\ref{sec:intro} each vertex-transitive graph of order $4p$
 has a $(4,p)$-semiregular automorphism. Using the program package {\sc Magma} \cite{Mag}
 a total of ten graphs of order $4p$ with a primitive automorphism group
 and having valency less  than $4p/3$, were found.
 For each of these graphs Table~\ref{tab:simbol} gives corresponding symbols
 by listing their entries $S_{i,j}$, $i,j \in \ZZ_4$.
 Among these graphs only the Coxeter graph is without a Hamilton cycle (the graph $X_2$ in Table~\ref{tab:simbol}). 
 This fact can be easily seen from the structure of the corresponding quotient graphs relative to
 a $(4,p)$-semiregular automorphism.
 Namely for each of these graphs, the quotient has a Hamilton cycle containing multiedges,  
 and so this cycle lifts to a Hamilton cycle in the  original graph (see also Proposition~\ref{pro:lem5}).

\begin{table}[h]
{\small
 $$
 \begin{array}{|c|c|c|c|c|c|c|c|c|c|c| }
\hline

& X_1&X_2&X_3&X_4&X_5&X_6&X_7&X_8&X_9&X_{10}\\
\hline
\hline
  |V(X_i)|& 28&28&28&28&28&28&68&68&68&52 \\
\hline
  valency& 9&3&6&6&9&9&12&15&20&6 \\
\hline
  S_{1,1}&\pm 3&\pm 1&\emptyset&\emptyset&\pm 2, \pm 3&\pm 1 &\pm 2,\pm 5 & \pm 6,\pm 8 &\pm 2,\pm 3, \pm 6  &\emptyset \\
\hline
  S_{2,2}& \emptyset&\pm 2& \pm 3&\emptyset&\pm 1, \pm 3&\pm 3 &\pm 1,\pm 6 &\pm  1, \pm 5 &\pm  5, \pm 7, \pm 8 &\emptyset \\
\hline
  S_{3,3}& \pm 1,\pm 3& \pm 4&\pm 1&\emptyset&\emptyset&\emptyset&\pm 4, \pm 7&\pm 2, \pm 7&\pm 4&\emptyset \\
\hline
  S_{4,4}&\pm 1 & \emptyset&\pm 2& \emptyset&\pm 1, \pm 2&\pm 2&\pm 3,\pm 8&\pm 3,\pm 4&\pm 1&\emptyset \\
\hline
  S_{1,2}&0,\pm 2 & \emptyset& 0,3 &0,3 &0 &0,4 &0,1 &0,5,7,9,14 &0,12,13,16  &0,4  \\
\hline
  S_{1,3}& 0,6& \emptyset&  0,6 & 0,6 &0,\pm 2 &0,1,4 &0,15 &0 & 0,\pm 1, 9,10,11 &0,10\\
\hline
  S_{1,4}& 0,4&0& 0,5 &0,5 &0 &0,1 &0,12,13,16 &0,1,2,6,13  &0,\pm 5, 10&0,12\\
\hline
  S_{2,3}&2,4 & \emptyset&5&0,3&0, \pm 3&0,2,4&2,4,10,12&6,\pm 7,13,14 &3,6,9,12&0,6\\
\hline
  S_{2,4}&0,1,\pm 3&0&1&0,2&0&1,3&1,10&11&\pm 2,4, 6,8,12&0,8\\
\hline
  S_{3,4}&1&0& 4&0,1& 0,\pm 1& 0,\pm 1& 8,12& 5,\pm 7, 8,16& \pm 2, 3, \pm 7, 9,12,13 & 0,11\\
\hline
\hline
 \end{array}
  $$
  \caption{\label{tab:simbol} {\footnotesize Symbols of connected vertex-transitive graphs of valency less than one third of
the number of vertices arising  from the actions in Proposition~\ref{pro:3.1}.}}}
\end{table}

For  the action of $A_8$ on cosets of $S_6$ and the action of
$S_8$ on cosets of $S_6\times \ZZ_2$ (part (i) of Proposition~\ref{pro:3.1})
the corresponding orbital graphs have valencies $12$ and $15$,  and thus more than $28/3$.
So these graphs are hamiltonian by Proposition~\ref{pro:jack}.

\begin{figure}[h!]
\begin{footnotesize}
\begin{center}
\includegraphics[width=0.45\hsize]{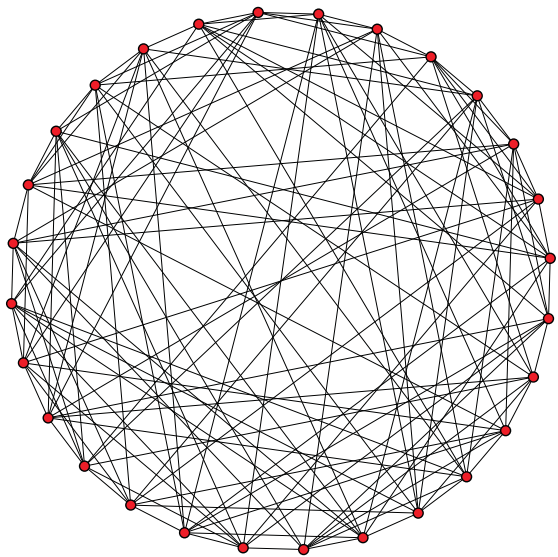}
\caption{\label{Prim28_2}\footnotesize The vertex-transitive
graph on $28$ vertices  with primitive automorphism group of
valency $9$, arising from the action of the group $PSL(2,8)$ on
the cosets of a subgroup $D_{18}$.}
\end{center}
\end{footnotesize}
\end{figure}

For the action of $PSL(2,8)$ on the cosets of $D_{18}$
(part (ii) of Proposition~\ref{pro:3.1})
we get that $D_{18}$  has three nontrivial suborbits, all of which are self-paired of length $9$.
Graphs arising from these suborbits are all isomorphic to the graph $X_1$ given in Table~\ref{tab:simbol}
(see also Figure~\ref{Prim28_2}).

For the action of $PGL(2,7)$ on the cosets of $D_{12}$ (part (iii) of Proposition~\ref{pro:3.1})
we deduce that $D_{12}$  has four nontrivial suborbits (all self-paired) 
one of which is of length $3$, two of length $6$ and one of length $12$. 
The graph arising from the  suborbit
of length $3$ is isomorphic to the Coxeter graph ($X_2$ in Table~\ref{tab:simbol}).
Next, $X_3$ and $X_4$ arise from the two suborbits of length $6$.
One of the graphs arising from the union of the suborbit of length $3$
and a suborbit of length $6$ is isomorphic to the graph $X_5$  and the other one to the graph $X_6$ in
Table~\ref{tab:simbol}.
As for the graph associated with  the suborbit
of length $12$,  it is clearly hamiltonian by  Proposition~\ref{pro:jack}.

The action of $PSL(2,16)\le G\le P\Gamma L(2,16)$
on cosets of $N_G(PGL(2,4))$ (part (iv) of Proposition~\ref{pro:3.1}), 
we deduce that $N_G(PGL(2,4))$  has four nontrivial suborbits, all
of which are self-paired, one of length $12$, one of length $15$ and two of length $20$.
The corresponding graphs  are, respectively, $X_7$, $X_8$ and $X_9$ in Table~\ref{tab:simbol}.

\begin{figure}[h]
\begin{footnotesize}
\begin{center}
\includegraphics[width=0.55\hsize]{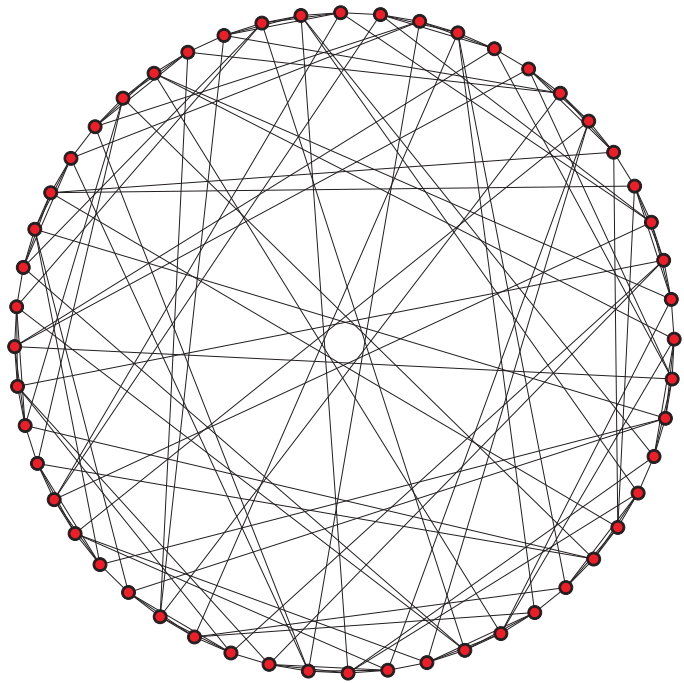}
\caption{\label{Prim52}\footnotesize The vertex-transitive
graph on $52$ vertices  with primitive automorphism group of
valency $6$, arising from the action of the group $\PSL(3,3)\le G\le \PSL(3,3)$
acting on the $52=4p$ incident point-line pairs of $PG(2,3)$.}
\end{center}
\end{footnotesize}
\end{figure}

As for the action of $PSL(3,3)\le G\le  PGL(3,3)$ on the $52$ incident point-line pairs of $PG(2,3)$
(part (v) of Proposition~\ref{pro:3.1}), we deduce that
there are five nontrivial suborbits, two of which are non-self-paired of length $9$, and three 
are self-paired of lengths $3$, $3$ and $27$.
The graph arising from the union of the two non-self-paired suborbits
has valency $18$ and is hamiltonian by  Proposition~\ref{pro:jack}, 
as is for the same reason the graph associated with the suborbit of length $27$.
The graphs arising from  the  suborbits of length $3$ are both disconnected.
Their union is isomorphic to the graph $X_{10}$ in Table~\ref{tab:simbol}
(see also Figure~\ref{Prim52}).

\end{document}